\documentclass{amsart}
\usepackage[cp1251]{inputenc}
\title{Local Lipschitz property of the map which puts in correspondence to $N$--net its Chebyshev centre.}
\author{P.~N.~Ivanshin}
\author{E.~N.~Sosov}
\address{Pyotr N. Ivanshin 
\newline\hphantom{iii}N.G. Chebotarev's RIMM 
\newline\hphantom{iii}Kazan State University.
\newline\hphantom{iii}420008, Kazan, Universitetskaya, 17}

\email{Pyotr.Ivanshin@ksu.ru}   

\address{Evgenii N. Sosov  
\newline\hphantom{iii}N.G. Chebotarev's RIMM
\newline\hphantom{iii}Kazan State University.
\newline\hphantom{iii}420008, Kazan, Universitetskaya, 17}

\email{Evgenii.Sosov@ksu.ru}   

\begin{document}
\maketitle

We prove local Lipschitz property of the map which puts in correspondence to each $N$--net different from $(N-1)$-net its Chebyshev center. If dimension of Euclidean or Lobachevsky space is greater than $1$ and net consists of more than $2$ points we show that this map is not Lipschits in a neighbourhood of the space of all $2$--nets embedded into the space of $N$--nets endowed by Hausdorff metric.

\vskip10pt
\section{Definitions and notation}
\vskip10pt

We assume the following notation. 
\newline 
$\mathbb{R}_{+}$ --- the set of all non-negative real numbers. 
\newline 
$(X,\rho)$ --- metric space with metric $\rho$.
\newline
$xy = \rho (x,y)$, $xZ = \inf \{xu : u \in Z \}$
for $x, \,y \in X$, $Z \subset X$.
\newline 
$\omega (x,y)$ --- the non-empty set such that for any $z \in \omega (x,y)$
\newline
$2 xz = 2 yz = xy$ for $x$, $y \in X$
(in the Euclidean space it is the middle-point of the interval $[x,y]$). 
\newline
$B[x,r]$ ($B(x,r)$, $S(x,r)$) --- a closed ball (open ball, sphere) with the center $x \in (X,\rho)$ and of radius $r \geq 0$.
\newline 
$\Sigma_N (X)$ --- the space of all non-empty subsets of $(X,\rho)$
consisting of no more than $N$ points. Elements of the space $\Sigma_N (X)$ are called  $N$--nets [1].
\newline 
$\alpha : \Sigma_N (X) \times  \Sigma_N (X) \rightarrow \mathbb{R}_{+}$, 
$\alpha  (M,T) = \max \{\max \{xT :  x \in M\},\max \{tM : t  \in T\}\}$~---
Hausdorff metric on the set $\Sigma_N (X)$ (cf. [2], p. 223). 
\newline
$B_{\alpha}(M,r)$ --- open ball centered in the point $M \in (\Sigma_N (X),\alpha)$, 
of the radius $r > 0$. 
\newline 
$D[M]$ --- diameter of the set $M \in \Sigma_N (X)$.
\newline  
Consider $M \in \Sigma_N (X)$. Assume that $R(M) = \inf\{\sup [yx : y \in M] : x \in X \}$ is Chebyshev radius of the $N$--net $M$. The point $\mathrm{cheb}(M) \in X$ is called Chebyshev center if 
$\sup [x\mathrm{cheb}(M) : x \in M] = R(M)$  [1].
\newline
We call a mapping $f : (X,\rho) \rightarrow (X_1,\rho_1)$ Lipschitz if there exists a constant $L \geq 0$ (it is called Lipschitz constant) such that for all $x$, $y \in X$
$\rho_1(f(x),f(y)) \leq \rho (x,y)$. We call Lipschitz map with the constant equal to $1$ non-expanding one (cf. [3], p. 10).

The notation given below will be used in Euclidean space as well as in Lobachevsky one.
\newline 
$co(M)$ --- convex hull of the set $M$.
\newline 
$[x,y]$ ($(x,y)$, $(x,y]$) --- closed (open, half-closed) interval of the end-points $x, y$. 
\newline
$\Pi (x_1,\ldots,x_{m+1})$ --- $m$--plane defined by points $x_1,\ldots,x_{m+1}$. 
\newline
$\Lambda (x,y)$ --- ray with the vertex in the point $x$ comprising $y \neq x$. 

\section{Statement of the problem}

The notion of Chebyshev center is fundamental thus there exists a vast choice of its application (approximation theory, functional analysis, geometry, technical sciences and such like). The authors study here only one local Lipschitz property of the mapping which puts in correspondence to each $N$-net (which is not a $N-1$--net) of Euclidean space its Chebyshev center.
This property despite its obvious importance was not completely considered even in the easiest case of Euclidean space. 

It is well-known that the Chebyshev center of the non-empty bounded subset of Euclidean or Lobachevsky space $X$ is strong continuous ([4], [5]), i.e. the mapping 
$\mathrm{cheb} : (B(X),\alpha) \rightarrow X$, $M \mapsto \mathrm{cheb}(M)$,
here $B(X)$ is the set of all non-empty bounded subsets of the space $X$,
is continuous. 
Moreover, restriction of the map $\mathrm{cheb}$ to the set of all balls of the space with inner metric is a non-expanding map [7]. Nevertheless the restriction of the aforementioned map $\mathrm{cheb}$ to the set of all non-empty bounded closed convex subsets of Euclidean plane is not Lipschitz even in a neighbourhood of a closed circle [8].

Thus we arrive to the following (rather complicated in general case) task:
Find classes of subsets (different from one consisting of the points or balls) of the given metric space such that for any class
$K$ the map $\mathrm{cheb} : (K,\alpha) \rightarrow X$, $M \mapsto \mathrm{cheb}(M)$ 
is locally Lipschitz.

It seems natural to solve this problem first for Euclidean space $(X,\rho)$ ($X$ may be infinite-dimensional) with the standard metric and  
$K = \Sigma_N (X)$ or $K = \Sigma_N (X)\backslash \Sigma_{N-1} (X)$.
The brief formulation of the result achieved by the authors is as follows:
We prove that the map 
$\mathrm{cheb} : (\Sigma_N (X)\backslash \Sigma_{N-1} (X),\alpha) \rightarrow X$, $M \mapsto \mathrm{cheb}(M)$ is locally Lipschitz for Euclidean space $X$. 
In case dimension of the Euclidean or Lobachevsky space is greater than $1$
and $N>2$ we show that this map is not Lipschitz in a neighbourhood of the subspace  $\Sigma_2 (X) \subset (\Sigma_N (X),\alpha)$. 

\section{Statement of the main results}

Here we proceed to the detailed formulation of our results.
If $N=2$ lemma 1 gives the answer to the problem stated in the previous paragraph. This lemma is a part of lemma 2 from [8].

{\bf Lemma 1.} Assume that for any points $p, \, x, \,y$ of the metric space   
$(X,\rho)$ the set $\omega (p,x)$ consists of one point. Assume also that the inequality
$$2 \omega (p,x)\omega (p,y) \leq  xy$$ holds true. 
Then the inequalities 
$$\mathrm{cheb}(M)\mathrm{cheb}(Z) \leq \alpha (M,Z) \leq \mathrm{cheb}(M)\mathrm{cheb}(Z) +
(D[M] + D[Z])/2$$ also hold 
for any $M,\, Z \in \Sigma_2 (X)$.

We give the answer to our question for Euclidean line $\mathbb{R}$ in lemma 2. 

{\bf Lemma 2.} Let $X$ be Euclidean line. Then the map
\newline
$\mathrm{cheb} : (\Sigma_N (X),\alpha) \rightarrow X$ is non-expanding.

The global Lipschitz property of the map $\mathrm{cheb}$ vanishes if the dimension of the space is greater than $1$ and $N > 2$.

{\bf Lemma 3.} Let $X$ be Euclidean space (Lobachevsky space) of dimension greater than $1$ and $N>2$. Let $U$ denote a neighbourhood of 
$\Sigma_2 (X) \subset \Sigma_N (X)$. 
Then
 
$(i)$ the map $\mathrm{cheb} : (U,\alpha) \rightarrow X$ is not Lipschitz;

$(ii)$ in Euclidean space the mapping $\mathrm{cheb} : (\Sigma_N (X),\alpha) \rightarrow X$
is not uniform continuous.

The following lemmas help us to prove local Lipschitz property of the map $\mathrm{cheb}$ in case $N = 3$.

{\bf Lemma 4.} Let $X$ be Euclidean space. Then
for any pair $M = \{u,v,w\}$ and $Z = \{u,v,z\} \in \Sigma_3 \backslash \Sigma_2$ such that $w \in (u,z)$ we get inequality 
$$
\mathrm{cheb}(M)\mathrm{cheb}(Z) \leq L \alpha (M,Z),
$$ 
here $L = 1/(2\sin (\varphi))$ if $\varphi = \angle vuw \neq 0$ and $\angle uvw$ are
acute angles; and $L = 1/2$ in other cases.

{\bf Lemma 5.}
Let  $X$ be Euclidean space, $M = \{u,v,w\} \in \Sigma_3 \backslash \Sigma_2$ and
$\varepsilon = \min [ab : a \neq b, a \, b \in M]/8$. Then there exists a constant $L > 0$ such that for all $z_1$, $z_2 \in B(w,\varepsilon)$
the following inequality holds true
$$
\mathrm{cheb}(W_1)\mathrm{cheb}(W_2) \leq L \alpha (W_1,W_2),
$$ 
here  
$W_1 =  \{u,v,z_1\}$, $W_2 = \{u,v,z_2\}$.

Now we are ready to prove local Lipschitz property of the map $\mathrm{cheb}$ in case $N = 3$.

{\bf Corollary 1.}
Let $X$ be Euclidean space, $M \in \Sigma_3 \backslash \Sigma_2$ and
$\varepsilon = \min [ab : a \neq b, a \, b \in M]/8$. Then 
$\mathrm{cheb} : (B_{\alpha}(M,\varepsilon),\alpha) \rightarrow X$ is  Lipschitz map.

Theorem 1 is the extension of this result to the case $N > 3$
in Euclidean plane. 

{\bf Theorem 1.} Let $X$ be Euclidean plane, $N > 3$,
$M \in \Sigma_N \backslash \Sigma_{N-1}$ and 
$\varepsilon = \min [ab : a \neq b, a \, b \in M]/8$. Then 
$\mathrm{cheb} : (B_{\alpha}(M,\varepsilon),\alpha) \rightarrow X$ is Lipschitz map.
  
Then we turn to the case of Euclidean space of dimension greater than $2$.


{\bf Lemma 6.} Assume that $X$ is a Euclidean space of dimension greater than $1$ and $N$--net $M = \{x_1,\ldots,x_N\} \in \Sigma_N \backslash \Sigma_{N-1}$ 
defines $(N-1)$--dimensional simplex. Then there exists a constant $L > 0$ such that 
for each $N$--net $Z = \{x_1,\ldots,x_{N-1},y_N\}$ meeting the condition $x_N \in (x_1,y_N)$ the inequality
$$
\mathrm{cheb}(M)\mathrm{cheb}(Z) \leq L \alpha (M,W)
$$ 
holds true.

{\bf Lemma 7} Let $X$ be Euclidean space of dimension greater than $1$,  $N$-net 
$M = \{x_1,\ldots,x_N\} \in \Sigma_N (X) \backslash \Sigma_{N-1} (X)$ define an $(N-1)$-dimensional simplex and $W = \{x_1,\ldots,x_N,x_{N+1}\} \in \Sigma_{N+1} (X) \backslash \Sigma_N (X)$ be such that
$x_{N+1} \in \Pi (x_1,\ldots,x_N)\backslash M$ and 
$$\varepsilon = \min [\min [ab : a \neq b, a, \, b \in \{x_1,\ldots,x_{N-1}\}],x_N \Pi (x_1,\ldots,x_{N-1})]/8.$$ 
Then 

$(i)$ there exists a constant $L > 0$ such that for all $z_1$, $z_2 \in B(x_N,\varepsilon)$
the inequality 
$$\mathrm{cheb}(Z_1)\mathrm{cheb}(Z_2) \leq L \alpha (Z_1,Z_2),
$$ holds true, here  
$Z_1 =  \{x_1,\ldots,x_{N-1},z_1\}$, $Z_2 = \{x_1,\ldots,x_{N-1},z_2\}$. 

$(ii)$ there exist constants $\delta>0$ and $L > 0$ such that for all $y_1$, $y_2 \in B(x_{N+1},\delta)$
the inequality 
$$
\mathrm{cheb}(Y_1)\mathrm{cheb}(Y_2) \leq L \alpha (Y_1,Y_2),
$$ 
holds, here  
$Y_1 =  \{x_1,\ldots,x_N,y_1\}$, $Y_2 = \{x_1,\ldots,x_N,y_2\}$. 

{\bf Theorem 2.} Let $(X,\rho)$  be Euclidean space.
Then for any $M \in \Sigma_N \backslash \Sigma_{N-1}$ there exists $\varepsilon = \varepsilon (M)$ such that the map 
$\mathrm{cheb} : (B_{\alpha}(M,\varepsilon),\alpha) \rightarrow X$ is a Lipschitz map.

It turns out that Lipschits property of the map $\mathrm{cheb}$ holds for two $N$--nets which lie sufficiently far from each other.

{\bf Statement 1.} Let $X$ be Euclidean plane and $N > 2$.
Then for any $M$, $Z \in \Sigma_N (X)$ such that 
$B(\mathrm{cheb}(M),R(M)) \cap B(\mathrm{cheb}(Z),R(Z)) = \emptyset$ the inequality
$$
\mathrm{cheb}(M)\mathrm{cheb}(Z) \leq L \alpha (M,Z)
$$ 
holds true,
here $L = (1 + \sqrt {5})/2$ for $N > 3$ and $L = \sqrt {2}$ for $N = 3$.

{\bf Statement 2.} 
$(i)$ Consider a Euclidean space $X$ and $3$--net $M = \{u,v,w\}$, 
\newline
Let $3$--net $Z = \{u,v,z\} \in \Sigma_3 (X) \backslash \Sigma_2 (X)$ be such that:

1.  $co(\{u,v,w\})\cap co(\{u,v,z\}) = [u,v]$;

2. if angles $\angle (uwv)$, $\angle (uzv)$ are acute then $\alpha (M,Z) < wz$.
\newline
Then the inequality 
$$
\mathrm{cheb}(M)\mathrm{cheb}(Z) \leq \alpha (M,Z)
$$ 
holds true.

$(ii)$ Let $X$ be Euclidean plane. Then the inequality 
$$
\mathrm{cheb}(M)\mathrm{cheb}(Z) \leq 2 \alpha (M,Z)
$$ 
holds for all 
$M = \{u,v,w\}$ and $Z = \{u,q,z\} \in \Sigma_3 (X) \backslash \Sigma_2 (X)$
such that
$$
co(\{u,v,w\})\cap co(\{u,v,z\}) = \{u\}.
$$   
\newline

\section{Proofs of the results}

{\bf Proof of lemma 2.}
Consider two arbitrary $N$--nets 
$M = \{x_1,\ldots,x_N\}$, $Z = \{y_1,\ldots,y_N\}$.
Let us  assume that $x_1 \leq \ldots \leq x_N$, $y_1 \leq \ldots \leq y_N$ and $x_1 \leq y_1$. 
Then lemma 1 together with the definitions of the Hausdorff metric and of the Chebyshev center gives us the inequality 
$$
\mathrm{cheb}(M)\mathrm{cheb}(Z) = \mathrm{cheb}(\{x_1,x_N\})\mathrm{cheb}(\{y_1,y_N\}) \leq \alpha (\{x_1,x_N\},\{y_1,y_N\})
\leq 
$$
$$
\leq \alpha (M,Z).
$$

{\bf Proof of lemma 3.} 
$(i)$ It suffices to verify the statement in case  $N=3$.
Let $x$, $y$ be two different points in Euclidean plane (Lobachevsky plane)
and $S_{+}$ be a semicircle constructed on the interval $[x,y]$ as diameter of the circle. Consider a non-fixed point $z \in S_{+}$ such that $yz \in (0,xy/2)$. Now find a point $u$ meeting the following conditions:
a) $z \in [x,u]$; b) the point $y$ lies on a circle based on $[x,u]$ as diameter. 
Consider $3$--nets 
$M = \{x, y, z\}$, $W = \{x, y, u\} \in \Sigma_3 \backslash \Sigma_2$ 
in Euclidean space (or in Lobachevsky space).
Then by construction $v = \mathrm{cheb} (M) = \omega (x,y)$, $w = \mathrm{cheb} (W) = \omega (x,u)$ and $uz = \alpha (M,W)$.
In Euclidean plane for any constant $L > 0$ we can choose a point $z \in S_{+}$ so that   $2L yz < xy$. 
Then
$\mathrm{cheb}(M)\mathrm{cheb}(W) = vw = (xy)(uz)/(2zy) > L uz = L \alpha (M,W)$.
This proves lemma 2 for Euclidean space.

Now let $p$ be a base of a perpendicular to the interval $[x,w]$ passing through point $v$ and $\psi$ be an angle in the vertex $w$ of the triangle $\triangle xvw$ in Lobachevsky plane. The triangles
$\triangle xvw$, $\triangle wpv$ are right-angled by construction. Hence,
we get formulae (cf. (4a), (4b) from [8], pg. 58)
$\tanh (pv) = \sinh (pw) \tan (\psi)$,
$\tanh (vx) = \sinh (vw) \tan (\psi)$.
Then the fraction   
$\sinh (vw)/\sinh (pw) = \tanh (vx)/\tanh (pv)$
tends to $\infty$ as ($pv \rightarrow 0$).
But then 
$\mathrm{cheb}(M)\mathrm{cheb}(W)/\alpha (M,W) =(vw)/(2pw)$ 
also infinitely increases as ($pv \rightarrow 0$).
This completes the proof for Lobachevsky space.
\newline
$(ii)$ Let us use the notation of the previous part of the proof. At the same time let us introduce new points $x_n = y + n (x - y)$, $z_n = \Lambda (u,x_n)\cap S(\omega (x_n,y), x_ny/2)$ and triples $M_n = \{x_n, y, z_n\}$, $Z_n = \{x_n, y, u\},$ here we identify points with their radius-vectors and $n =1,2,\ldots$. Then $\alpha (M_n,Z_n) \rightarrow 0$ but $\mathrm{cheb}(M_n)\mathrm{cheb}(Z_n) \rightarrow yu$
as ($n \rightarrow \infty$).
This completes the proof of lemma 3.

{\bf Proof of lemma 4.}
Let $\varphi$ be a nonzero acute angle.
Then
$$
\mathrm{cheb}(M)\mathrm{cheb}(W) \leq wz/2 = alpha (M,W)/2.
$$
Let $\varphi \neq 0$ be acute angle, 
$p$ --- base of the perpendicular to the ray $\Lambda (u,w)$ passing through $v$, $q$ --- point on the ray $\Lambda (u,w)$ such that the interval $[q,v]$ is perpendicular to the interval $[u,v]$.
Then the set of points of the ray $\Lambda (u,w)$ can be represented as follows: 
$\Lambda (u,w) = (u,p]\cup (p,q)\cup (\Lambda (u,w) \backslash (u,q))$.
Let us consider position of $\mathrm{cheb}(M)$ depending on $w$. 
If $w \in (u,p]$ then $c = \mathrm{cheb}(M) = \omega (u,v)$. If $w \in (p,q)$ then 
$\mathrm{cheb}(M) \in (c,b)$, here the point $b \in (u,q]$ is such that the interval $[b,c]$
is perpendicular to the ray $\Lambda (u,v)$ and $\mathrm{cheb}(M) = pw/(2\sin (\varphi))$.
If $w \in \Lambda (u,w) \backslash (u,q)$ then $\mathrm{cheb}(M) = \omega (u,w)$.
Let us investigate all possibilities for location of the points $w$, $z$ on the ray $\Lambda (u,w)$.

If either $w$, $z \in (u,p]$ or $w$, $z \in \Lambda (u,w) \backslash (u,q)$ then
$\mathrm{cheb}(M)\mathrm{cheb}(W) = 0$.
 
In case $w \in (u,p)$, $z \in (p,q)$ we have
$\mathrm{cheb}(M)\mathrm{cheb}(W) = pz/(2\sin (\varphi)) = \alpha (M,W)/(2\sin (\varphi))$.

If $w \in (u,p)$, $z \in \Lambda (u,w) \backslash (u,q)$ then
$\mathrm{cheb}(M)\mathrm{cheb}(W) = c\mathrm{cheb}(W) \leq cb + b\mathrm{cheb}(W) = pq/(2\sin (\varphi)) + qz/2 
\leq pz/(2\sin (\varphi)) \leq \alpha (M,W)/(2\sin (\varphi))$.

If $w$, $z \in (p,q)$ then
$\mathrm{cheb}(M)\mathrm{cheb}(W) = (pz - pw)/(2\sin (\varphi)) \leq \alpha (M,W)/(2\sin (\varphi))$.

If $w \in (p,q)$, $z \in \Lambda (u,w) \backslash (u,q)$ then
$\mathrm{cheb}(M)\mathrm{cheb}(W) \leq \mathrm{cheb}(M)b + b\mathrm{cheb}(W) = (pq - pw)/(2\sin (\varphi)) + qz/2 = 
\alpha (M,W)/(2\sin (\varphi))$.

Thus in all considered cases we can fix $1/(2\sin (\varphi))$ as Lipschitz constant $L$. This completes the proof of the lemma.

{\bf Proof of lemma 5.}
Let triangles defined by $3$--nets $W_1$,  $W_2$ be not acute. Then with the help of lemma 1 we get the desired inequality with constant $L = 1/2$.  Now let $X$ be a Euclidean plane and 
$3$--net $W_1$ define an acute triangle. If the angle 
$\angle uz_2v$ is blunt then there exists a point $z_3 \in (z_1,z_2)$ such that the angle  $\angle uz_3v$ is right and $\mathrm{cheb} (\{u,v,z\}) = \mathrm{cheb} (W_2)$. Thus without loss of generality one can assume that the triangle given by the $3$--net $W_2$ is not blunt. 
Let us consider all possible cases of the location of the point $z_2 \in B(w,\varepsilon)$.

1. Let a point $z_2\neq z_1$ belong to the set
$co \{u,v,z_1\}$ or to the closed angle with vertex $z_1$ and sides on the rays 
$\Lambda (v,z_1)$, $\Lambda (u,z_1)$. If $z_2 \in \Lambda (u,z_1)\cup \Lambda (v,z_1)$ then we get the desired inequality from lemma 4. In the other case let us denote by $p$ a point of intersection of the rays $\Lambda (v,z_2)$, 
$\Lambda (u,z_1)$ and $W_3 = \{u,v,p\}$. Note that the angle $\angle z_1pz_2$ 
is not acute. 
Now we use lemma 4 and triangle inequality. Thus we get constants $L_1 = \sup \{1/(2\sin (\angle z_1uv)) : z_1 \in B(w,\varepsilon)\} < \infty$ and
$L_2 = \sup \{1/(2\sin (\angle z_2vu)) : z_2 \in B(w,\varepsilon)\} < \infty$ such that   
$$
\mathrm{cheb}(W_1)\mathrm{cheb}(W_2) \leq \mathrm{cheb}(W_1)\mathrm{cheb}(W_3) + \mathrm{cheb}(W_3)\mathrm{cheb}(W_2) \leq
L_1 \alpha (W_1,W_3) +
$$
$$
+ L_2 \alpha (W_3,W_2) \leq
(L_1  + L_2) \alpha (W_1,W_2).
$$
So, in this case $L = L_1  + L_2$.

2. Now let a point $z_2\neq z_1$ belong to the open angle with vertex $z_1$ and sides on the rays $\Lambda (z_1,u)$ and $\Lambda (v,z_1)$ (the case in which $z_2\neq z_1$ belongs to the open acute angle with vertex $z_1$ and sides on the rays $\Lambda (z_1,v)$ and $\Lambda (u,z_1)$ is similar to this one). Let us denote by $p$ the intersection point of the rays $\Lambda (v,z_2)$, $\Lambda (u,z_1)$ and $W_3 = \{u,v,p\}$.  If the angle $\angle z_1pz_2$ is not acute we can use the same considerations as in the first case. 
If the angle $\angle z_1pz_2$ is acute then we have inequalities
$\alpha (W_1,W_3) = pz_1 = pv \sin (\angle pvz_1)/\sin (\angle uz_1v) 
\leq vz_1 \sin (\angle pvz_1)/\sin (\angle uz_1v) \leq z_1z_2/\sin (\angle uz_1v) 
\leq L_3  \alpha (W_1,W_2)$,
here $L_3 = \sup \{1/(\sin (\angle uz_1v)) : z_1 \in B(w,\varepsilon)\} < \infty$ by lemma 5.   
Now we again use this inequality, triangle inequality and lemma 4. Thus we get
$$
\mathrm{cheb}(W_1)\mathrm{cheb}(W_2) \leq \mathrm{cheb}(W_1)\mathrm{cheb}(W_3) + \mathrm{cheb}(W_3)\mathrm{cheb}(W_2)
\leq L_1 \alpha (W_1,W_3) +
$$
$$
+ L_2 \alpha (W_3,W_2) \leq 
L_1 \alpha (W_1,W_3) + L_2 (\alpha (W_1,W_3) + \alpha (W_1,W_2)) \leq L \alpha (W_1,W_2),
$$
here $L = L_3 (L_1 + L_2) + L_2$, 
$L_1 = \sup \{1/(2\sin (\angle z_1uv)) : z_1 \in B(w,\varepsilon)\} < \infty$, 
$L_2 = \sup \{1/(2\sin (\angle z_2vu)) : z_2 \in B(w,\varepsilon)\} < \infty$. 

Now let $X$ be Euclidean space of dimension greater than $2$ and 
$W_4$ be an orthogonal projection of $W_2$ on the plane containing $W_1$; note that if $W_1$ lies on the line then $W_4 = W_2$. Then the previous considerations provide us with the constant $L_4 > 0$ such that for all 
$z_1$, $z_2 \in B(w,\varepsilon)$
the inequality 
$$
\mathrm{cheb}(W_1)\mathrm{cheb}(W_4) \leq L_4 \alpha (W_1,W_4) \leq  L_4 \alpha (W_1,W_2)
$$
holds true.
Note that the inequality 
$$
\mathrm{cheb}(W_4)\mathrm{cheb}(W_2) \leq \alpha (W_1,W_2)
$$ 
follows from purely geometric considerations. It suffices now to apply triangle inequality to complete the proof. 

{\bf Proof of corollary 1.}
Fix an arbitrary pair of nets $W_1 =  \{x,y,z\}$, $W_2 = \{u,v,w\} \in B_{\alpha}(M,\varepsilon)$ and put
$W_3 =  \{x,y,w\}$, $W_4 = \{x,v,w\}$. Then lemma 5, definition of the Hausdorff metric and triangle inequality provide us with the constants 
$L_1 \, L_2 \, L_3 >0$ such that
$$
\mathrm{cheb}(W_1)\mathrm{cheb}(W_2) \leq \mathrm{cheb}(W_1)\mathrm{cheb}(W_3) + \mathrm{cheb}(W_3)\mathrm{cheb}(W_4) +
$$
$$
+
\mathrm{cheb}(W_4)\mathrm{cheb}(W_2) \leq L_1 \alpha (W_1,W_3) + L_2 \alpha (W_3,W_4) +
L_3 \alpha (W_4,W_2) \leq
$$
$$
\leq (L_1 + L_2 + L_3) \alpha (W_1,W_2).
$$
Thus 
$\mathrm{cheb} : (B_{\alpha}(M,\varepsilon),\alpha) \rightarrow X$ is a Lipschitz map. This completes the proof.

{\bf Proof of theorem 1.}

$(i)$ Let us consider two arbitrary $N$--nets of the special kind 
$Z_1 = \{x_1,x_2,\ldots,x_N\}$, $Z_2 = \{y_1,x_2,\ldots,x_N\} \in B_{\alpha}(M,\varepsilon)$ and
introduce the parametrisation $x = x(s)$ by the length of the interval $[x_1,y_1]$ so that $x_1 = x(0)$, $y_1 = x(x_1y_1)$.
If $z \in (x_1,y_1]\cap B[\mathrm{cheb}(Z_1),R(Z_1)]$ then
$\mathrm{cheb}(Z_1) = \mathrm{cheb} (\{z,x_2,\ldots,x_N\}$. Hence 
we may assume that $x_1 \in S(\mathrm{cheb}(Z_1),R(Z_1))$.
Let us put in correspondence to any $s \in [0,x_1y_1]$ an $N$--net
$Z(s) = \{x(s),x_2,\ldots,x_N\}$ and a convex polygon $Q(s)$, 
given by vertices $Z(s)\cap S(\mathrm{cheb}(M(s)),R(M(s)))$.
For any $s \in [x_1y_1]$ divide the polygon $Q(s)$ into triangles  $\{\bigtriangleup_1(s),\ldots,\bigtriangleup_{k(s)}(s)\}$
with common vertex $x(s)$ by connecting this vertex by intervals with all other vertices of $Q(s)$. Now for each 
$s \in [x_1y_1]$ $\mathrm{cheb}(Q(s)) = \mathrm{cheb}(M(s))$ belongs either to the interior of one of the acute triangles or to the side common to two adjacent triangles (it is possible that these triangles are one and the same or even they are just points on an interval)
from the division $\{\bigtriangleup_1(s),\ldots,\bigtriangleup_{k(s)}(s)\}$.
Consider two cases.

1. Let $\mathrm{cheb}(Z_1)$ belong to the interior of the triangle 
$\bigtriangleup x_1ab \in \{\bigtriangleup_1(0),\ldots,\bigtriangleup_{k(0)}(0)\}$.
Then by the strong continuity of the Chebyshev center one can find a minimal number $s_1 \in (0,x_1y_1]$ such that either
$\mathrm{cheb}(Z_2) = \mathrm{cheb}(\{y_1,a,b\})$ in case $s_1 = x_1y_1$,
or $\mathrm{cheb}(Z(s_1)) = \mathrm{cheb}(\{x(s_1),a,b\})$ lies on the side of one or two of triangles of the division 
$\{\bigtriangleup_1(s_1),\ldots,\bigtriangleup_{k(s_1)}(s_1)\}$ in case 
$s_1 \in (0,x_1y_1)$.

2. Let $\mathrm{cheb}(Z_1)$ belong to the common side of the triangles
$\bigtriangleup x_1ab$, $\bigtriangleup x_1ac \in \{\bigtriangleup_1(0),\ldots,\bigtriangleup_{k(0)}(0)\}$.
Then the strong continuity of the Chebyshev center implies the existence of a minimal $s_1 \in (0,x_1y_1]$ such that either
$\mathrm{cheb}(Z_2) = \mathrm{cheb}(\{y_1,a,b\})$ or $\mathrm{cheb}(Z_2) = \mathrm{cheb}(\{y_1,a,c\})$ 
if $s_1 = x_1y_1$,
or $\mathrm{cheb}(Z(s_1)) = \mathrm{cheb}(\{x(s_1),a,b\})$ or
$\mathrm{cheb}(Z(s_1)) = \mathrm{cheb}(\{x(s_1),a,c\})$ lies on the side of one or two of triangles of the division $\{\bigtriangleup_1(s_1),\ldots,\bigtriangleup_{k(s_1)}(s_1)\}$ in case
$s_1 \in (0,x_1y_1)$.
 
In either of the considered cases lemma 5 and compactness of the interval $[x_1, y_1]$ imply the existence of the constant $L_1 >0$ such that
$\mathrm{cheb}(Z_1)\mathrm{cheb}(Z(s_1)) \leq L_1 \alpha (Z_1,Z(s_1))$.
If $s_1 \neq x_1y_1$ then for $\mathrm{cheb}(Z(s_1))$ we can proceed in similar way and get $s_2$. Continuing this process we get from lemma 5 constants 
$0 < s_1 < \ldots < s_m = x_1y_1$, $L_1,\ldots,L_m > 0$ such that  
$$
\mathrm{cheb}(Z_1)\mathrm{cheb}(Z_2) \leq \mathrm{cheb}(Z_1)\mathrm{cheb}(Z(s_1)) + 
\mathrm{cheb}(Z(s_1))\mathrm{cheb}(Z(s_2)) + \ldots +
$$
$$
+ \mathrm{cheb}(Z(s_{m-1}))\mathrm{cheb}(Z_2) \leq 
L_1 \alpha (Z_1,Z(s_1)) + \ldots + L_m \alpha (Z(s_{m-1}),Z_2) 
\leq
$$
$$
\leq (L_1 + \ldots + L_m) \alpha (Z_1,Z_2).
$$
This completes the proof of the lemma 6 in this special case.
  
$(ii)$ Choose arbitrary nets $Z_1 = \{x_1,x_2,\ldots,x_N\}$, $Z_2 = \{y_1,y_2,\ldots,y_N\} \in B_{\alpha}(M,\varepsilon)$ and fix $Z_3 =  \{y_1,x_2,\ldots,x_N\},\ldots,Z_{N+1} = \{y_1,y_2,\ldots,y_{N-1},x_N\}$. Then triangle inequality, lemma 5 and the definition of the Hausdorff metric give us constants $L_1 \,\ldots, L_{N-1} >0$ such that
$$
\mathrm{cheb}(Z_1)\mathrm{cheb}(Z_2) \leq \mathrm{cheb}(Z_1)\mathrm{cheb}(Z_3) + \mathrm{cheb}(Z_3)\mathrm{cheb}(Z_4) + \ldots
+ 
$$
$$
+\mathrm{cheb}(Z_{N+1})\mathrm{cheb}(Z_2) \leq L_1 \alpha (Z_1,Z_3) + L_2 \alpha (Z_3,Z_4) + \ldots
+ 
$$
$$
+ L_{N-1} \alpha (Z_{N+1},Z_2) \leq (L_1 + \ldots + L_{N-1}) \alpha (Z_1,Z_2).
$$
Thus 
$\mathrm{cheb} : (B_{\alpha}(M,\varepsilon),\alpha) \rightarrow X$ is a Lipschitz map.

{\bf Proof of lemma 6.}
The proof is by induction.
Lemma 4 implies consistence of the statement in case $N = 3$. Assume now that lemma 6 holds true for all numbers less or equal than $N-1$.
Let us show that the statement holds true for $N$.
Let $p \neq x_1$ be an intersection point of the ray $\Lambda (x_1,x_N)$
with the sphere $S(\mathrm{cheb}(\{x_1,\ldots,x_{N-1}\}),R(\{x_1,\ldots,x_{N-1}\}))$.
and $q$ be a point on $\Lambda (x_1,x_N)$ such that there exists $k \in \{2,\ldots,N\}$ such that 
$$
x_k \in S(\mathrm{cheb}(\{x_1,\ldots,\hat x_k,\ldots,x_{N-1},q\}),R(\{x_1,\ldots,\hat x_k,\ldots,x_{N-1},q\})),
$$
(here cap on the symbol means that this symbol must be excluded) closest to $x_1$. 
Then the set of points of the ray $\Lambda (x_1,x_N)$ can be represented as follows:
$\Lambda (x_1,x_N) = (x_1,p]\cup (p,q)\cup (\Lambda (x_1,x_N) \backslash (x_1,q))$.
Consider all possible locations of the points $x_N$, $y_N$ on $\Lambda (x_1,x_N)$.

If $x_N$, $y_N \in (x_1,p]$ then
$\mathrm{cheb}(M)\mathrm{cheb}(Z) = 0$.
If $x_N$, $y_N \in (\Lambda (x_1,x_N) \backslash (x_1,q))$ then there exist $k$, $j \in \{2,\ldots,N\}$ such that
$\mathrm{cheb}(M) =\mathrm{cheb}(\{x_1,\ldots,\hat x_k,\ldots,x_{N}\})$, 
$\mathrm{cheb}(Z) =\mathrm{cheb}(\{x_1,\ldots,\hat x_j,\ldots,x_{N-1},y_{N}\})$.
Moreover by induction assumption there exists a constant $L > 0$ such that
$$
\mathrm{cheb}(M)\mathrm{cheb}(Z) \leq 
L \alpha (\{x_1,\ldots,\hat x_k,\ldots,x_{N}\},\{x_1,\ldots,\hat x_j,\ldots,x_{N-1},y_{N}\})
\leq
$$
$$
\leq L \alpha(M,Z).
$$

If $x_N \in (x_1,q)$, $y_N \in (p,q)$ then
$\mathrm{cheb}(M)\mathrm{cheb}(Z) \leq \alpha (M,Z)/(2\cos (\psi))$,
here $\psi$ is an angle between the ray $\Lambda (x_1,x_N)$ and the normal to the scale $(x_1,\ldots,x_{N-1})$ of the simplex. 

Let $b$ be an intersection-point of the boundary of the simplex $(x_1,\ldots,x_{N})$
with normal to the scale $(x_1,\ldots,x_{N-1})$ passing through the point $\mathrm{cheb}(\{x_1,\ldots,x_{N-1}\})$ which does not belong to this face.
If $x_N \in (x_1,q)$, $y_N \in (\Lambda (x_1,x_N) \backslash (x_1,q))$ then there exist $k$, $j \in \{2,\ldots, N\}$ such that
$b =\mathrm{cheb}(\{x_1,\ldots,\hat x_k,\ldots,x_{N-1},q\})$, 
$\mathrm{cheb}(Z) =\mathrm{cheb}(\{x_1,\ldots,\hat x_j,\ldots,x_{N-1},y_{N}\})$.
Moreover by assumption there exist a constant $L_1 > 0$ such that
$$
\mathrm{cheb}(M)\mathrm{cheb}(Z) \leq \mathrm{cheb}(M) b + b\mathrm{cheb}(Z) \leq 
q\{p,x_N\}/(2\cos \psi) + 
$$
$$
+\mathrm{cheb}(\{x_1,\ldots,\hat x_k,\ldots,x_{N-1},q\})\mathrm{cheb}(\{x_1,\ldots,\hat x_j,\ldots,x_{N-1},y_{N}\}) 
\leq
$$
$$
\leq
q\{p,x_N\}/(2\cos \psi) +
$$
$$
+ 
L_1 \alpha (\{x_1,\ldots,\hat x_k,\ldots,x_{N-1},q\},\{x_1,\ldots,\hat x_j,\ldots,x_{N-1},y_{N}\})
\leq L \alpha(M,Z),
$$
here $L = (1 + L_1)/(2\cos \psi)$. This completes the proof

{\bf Proof of lemma 7.}
The proof is by induction.
 
$(i)$ Lemma 5 implies that statement $(i)$ of lemma 7 holds true for $N = 3$. Assume now that it holds for all numbers less or equal than $N - 1$ and prove it for $N$.
Consider two cases.
 
1. Let there exist $k \in \{1,\ldots,N-1\}$ such that
$pz_1 \leq z_1z_2 $, here $p$ is an intersection point of the ray $\Lambda (x_k,z_2)$ with the $(N-2)$-plane 
$\Pi (\{x_1,\ldots,\hat x_k,\ldots,x_{N-1},z_1\})$.
 
Let us denote by $Z_3$ the set $\{x_1,\ldots,x_{N-1},p\}$. Now lemma 6 implies the existence of a constant 
$L_2 > 0$ such that
$\mathrm{cheb}(Z_3)\mathrm{cheb}(Z_2) \leq L_2 \alpha (Z_3,Z_2)$.
 
By induction assumption there exists a constant  
$L_1 > 0$ such that
$$
\mathrm{cheb}(\{x_1,\ldots,\hat x_k,\ldots,x_{N-1},z_1\})\mathrm{cheb}(\{x_1,\ldots,\hat x_k,\ldots,x_{N-1},p\}) \leq
$$
$$
L_1 \alpha (\{x_1,\ldots,\hat x_k,\ldots,x_{N-1},z_1\},\{x_1,\ldots,\hat x_k,\ldots,x_{N-1},p\})
\leq L_1 \alpha (Z_1,Z_3).
$$
Let $\psi$ stand for an angle between normal vector to the scale $(x_1,\ldots,x_{N-1})$ directed into the simplex and normal vector to the scale $(x_1,\ldots,\hat x_k,\ldots,x_{N-1},z_1)$ directed outside of the simplex. Now purely geometrical considerations imply the inequality 
$$
\mathrm{cheb}(Z_1)\mathrm{cheb}(Z_3) \leq
$$
$$
\leq \max [1/\sin {\psi},1] \mathrm{cheb}(\{x_1,\ldots,\hat x_k,\ldots,x_{N-1},z_1\})\mathrm{cheb}(\{x_1,\ldots,\hat x_k,\ldots,x_{N-1},p\}).
$$
  
Thus 
$$
\mathrm{cheb}(Z_1)\mathrm{cheb}(Z_2) \leq \mathrm{cheb}(Z_1)\mathrm{cheb}(Z_3) + 
\mathrm{cheb}(Z_3)\mathrm{cheb}(Z_2) \leq
$$
$$
\leq L_1 \max [1/\sin {\psi},1] \alpha (Z_1,Z_3) + L_2 \alpha (Z_3,Z_2) \leq
(L_1 \max [1/\sin {\psi},1]  +  2 L_2) \alpha (Z_1,Z_2).
$$
2. Let the assumption of case 1 be wrong, $j \in \{1,\ldots,N-1\}$ and denote by $p_j$ the intersection point of the ray $\Lambda (x_j,z_2)$ and the $(N-2)$-plane $\Pi (x_1,\ldots,\hat x_j,\ldots,x_{N-1},z_1)$.
Now choose $k \in \{1,\ldots,N-1\}$ such that 
$p_kz_1 = \min [p_jz_1 : j \in \{1,\ldots,N-1\}]$ and consider
$Z_3 = \{x_1,\ldots,x_{N-1},p_k\}$. 
 
Note that the angle $\angle z_1p_kz_2$ is acute and introduce the inequality 
$$
\alpha (Z_1,Z_3) = p_kz_1 = p_kx_k \sin {\angle p_kx_kz_1}/\sin {\angle p_kz_1x_k} 
\leq
$$
$$
\leq x_kz_1 \sin {\angle p_kx_kz_1}/\sin {\angle p_kz_1x_k} \leq z_1z_2/\sin {\angle p_kz_1x_k} 
\leq L_3  \alpha (Z_1,Z_2),
$$
here general assumptions imply that $L_3 = \sup \{1/\sin {\angle p_kz_1x_k} : z_1 \in B(x_N,\varepsilon)\} < \infty$.   
This inequality, triangle one, lemma 6 and similar estimate for $\mathrm{cheb}(Z_1)\mathrm{cheb}(Z_3)$ of the first case put together give us constants $L_1$, $L_2 > 0$ such that 
$$
\mathrm{cheb}(Z_1)\mathrm{cheb}(Z_2) \leq \mathrm{cheb}(Z_1)\mathrm{cheb}(Z_3) + \mathrm{cheb}(Z_3)\mathrm{cheb}(Z_2)
\leq
$$
$$
\leq L_1 \alpha (Z_1,Z_3) + L_2 \alpha (Z_3,Z_2) \leq 
L_1 \alpha (Z_1,Z_3) + L_2 (\alpha (Z_1,Z_3) + \alpha (Z_1,Z_2)) \leq L \alpha (Z_1,Z_2),
$$
here $L = L_3 (L_1 + L_2) + L_2$.
Thus statement $(i)$ of lemma 7 holds true.
 
$(ii)$ Lemma 5 provides us with the proof of statement $(ii)$ of lemma 7 for $N = 2$. Let us assume that statement $(ii)$ holds true also for all natural numbers less or equal than $N - 1$ and prove it for $N$. Let $0 < \delta < \min [ab : a \neq b, a, \, b \in W]/8$.
Consider three cases.
 
1. Let $[y_1,y_2]\cap\Pi (x_1,\ldots,x_N) = \oslash$. Then our statement holds by compactness of the interval $[y_1,y_2]$ and previously proved statement $(i)$ of lemma 7.
 
2. Let $[y_1,y_2] \subset \Pi (x_1,\ldots,x_N)$. Introduce the natural parametrisation $y = y(s)$  of the interval $[y_1,y_2]$ by length so that $y_1 = y(0)$, $y_2 = y(y_1y_2)$.
If $y \in (y_1,y_2]\cap B[\mathrm{cheb}(Y_1),R(Y_1)]$ then
$\mathrm{cheb}(Y_1) = \mathrm{cheb} (\{x_1,x_2,\ldots,x_N,y\})$. 
 
Hence we can assume that $y_1 \in S(\mathrm{cheb}(Y_1),R(Y_1))$. Consider for any  $s \in [0,y_1y_2]$  $(N+1)$-net
$Y(s) = \{x_1,x_2,\ldots,x_N,y(s)\}$ and convex polygon  $Q(s)$ with vertices $Y(s)\cap S(\mathrm{cheb}(Y(s)),R(Y(s)))$. Now we divide polygon $Q(s)$ for each $s \in [x_1,y_1]$ into simplices $\{\bigtriangleup_1(s),\bigtriangleup_{k(s)}(s)\}$
with the common vertex $y(s)$, here $1\leq k(s)\leq 2$. For any $s \in [y_1,y_2]$ $\mathrm{cheb}(Q(s)) = \mathrm{cheb}(Y(s))$ belongs either to the interior of one simplex or to the common face of two simplices of the division $\{\bigtriangleup_1(s),\bigtriangleup_{k(s)}(s)\}$ (the degenerate cases are also possible).

Consider two cases.

A. Let $\mathrm{cheb}(Y_1)$ belong to the interior of the simplex
  
$(a_1,\ldots,a_{N-1},y_1) \in \{\bigtriangleup_1(0),\bigtriangleup_{k(0)}(0)\}$.
 
Then---since the Chebyshev center depends on the set is strong continuous way---one can find a minimal $s_1 \in (0,y_1y_2]$ such that either
$\mathrm{cheb}(Y_2) = \mathrm{cheb}(\{a_1,a_2,\ldots,a_{N-1},y_2\})$ for $s_1 = y_1y_2$
or 
$\mathrm{cheb}(Y(s_1)) = \mathrm{cheb}(\{y(s_1),x_2,\ldots,x_N\})$ belongs to the scale of one or two simplices of the division 
$\{\bigtriangleup_1(s_1),\bigtriangleup_{k(s_1)}(s_1)\}$ for 
$s_1 \in (0,y_1y_2)$.

B. Let $\mathrm{cheb}(Y_1)$ belong to the common face of the simplices 
$(a_1,\ldots,a_{N-2},b,y_1)$, 
$(a_1,\ldots,a_{N-2},c,y_1) \in \{\bigtriangleup_1(0),\bigtriangleup_{k(0)}(0)\}$.
 
Then again strong continuity of the Chebyshev center implies the existence of the minimal $s_1 \in (0,y_1y_2]$, such that either
$\mathrm{cheb}(Y_2) = \mathrm{cheb}(\{a_1,\ldots,a_{N-2},b,y_2\})$ or 
$\mathrm{cheb}(Y_2) = \mathrm{cheb}(\{a_1,\ldots,a_{N-2},c,y_2\})$ 
for $s_1 = y_1y_2$
or
$\mathrm{cheb}(Y(s_1)) = \mathrm{cheb}(\{a_1,\ldots,a_{N-2},b,y(s_1)\})$ or 
$\mathrm{cheb}(Y(s_1)) = \mathrm{cheb}(\{a_1,\ldots,a_{N-2},c,y(s_1)\})$ belongs to the scale of one or two simplices of the division
$\{\bigtriangleup_1(s_1),\bigtriangleup_{k(s_1)}(s_1)\}$ for
$s_1 \in (0,y_1y_2)$.

In either case compactness of the interval $[y_1, y_2]$ and statement $(i)$ of lemma 7 or induction hypothesis provide us with the constant $L_1 >0$ such that 
$\mathrm{cheb}(Y_1)\mathrm{cheb}(Y(s_1)) \leq L_1 \alpha (Y_1,Y(s_1))$.
If $s_1 \neq y_1y_2$ then similar considerations applied to  $\mathrm{cheb}(Y(s_1))$ give us $s_2$. Continuing this process we get the set of constants 
$0 < s_1 < \ldots < s_i = y_1y_2$, $L_1,\ldots,L_i > 0$ such that 
$$
\mathrm{cheb}(Y_1)\mathrm{cheb}(Y_2) \leq \mathrm{cheb}(Y_1)\mathrm{cheb}(Y(s_1)) + 
\mathrm{cheb}(Y(s_1))\mathrm{cheb}(Y(s_2)) + \ldots +
$$
$$
+\mathrm{cheb}(Y(s_{i-1}))\mathrm{cheb}(Y_2) \leq 
L_1 \alpha (Y_1,Y(s_1)) + \ldots + L_i \alpha (Y(s_{i-1}),Y_2) 
\leq 
$$
$$
\leq (L_1 + \ldots + L_i) \alpha (Y_1,Y_2).
$$
Thus our statement holds true in this special case.
 
3. There exists a unique point $u$ such that $u = [y_1,y_2]\cap\Pi (x_1,\ldots,x_N)$.
 
Since triangle inequality holds true it suffices to prove the statement for the interval $[y_1,u]$. Let $\tilde y_1$ denote orthogonal projection of $y_1$ onto $\Pi (x_1,\ldots,x_N)$. Then geometrical considerations give us the inequality
$$
\mathrm{cheb}(Y_1)\mathrm{cheb} (\{x_1,x_2,\ldots,x_N,\tilde y_1\}) \leq \alpha (Y_1, \{x_1,x_2,\ldots,x_N,\tilde y_1\}).
$$
Now case 2 of the proof implies existence of the constant $L >0$ such that 
$$
\mathrm{cheb}(Y_1)\mathrm{cheb}(Y_2) \leq \mathrm{cheb}(Y_1)\mathrm{cheb}(\{x_1,x_2,\ldots,x_N,\tilde y_1\}) +
$$
$$ 
+\mathrm{cheb}(\{x_1,x_2,\ldots,x_N,\tilde y_1\})\mathrm{cheb}(\{x_1,x_2,\ldots,x_N,u\}) \leq \alpha (Y_1, \{x_1,x_2,\ldots,x_N,\tilde y_1\}) +
$$
$$ 
+ L \alpha (\{x_1,x_2,\ldots,x_N,\tilde y_1\},\{x_1,x_2,\ldots,x_N,u\}) \leq (1 + L)\alpha (Y_1, \{x_1,x_2,\ldots,x_N,u\}).
$$
This completes the proof of the lemma.

{\bf Proof of theorem 2.}

$(i)$ Let us consider two arbitrary $N$--nets of the special kind: 
$Z_1 = \{x_1,x_2,\ldots,x_N\}$, $Z_2 = \{y_1,x_2,\ldots,x_N\} \in B_{\alpha}(M,\varepsilon)$ and introduce a parametrisation $x = x(s)$ by the length of the interval 
$[x_1,y_1]$ such that $x_1 = x(0)$, $y_1 = x(x_1y_1)$.
If $z \in (x_1,y_1]\cap B[\mathrm{cheb}(Z_1),R(Z_1)]$ then
$\mathrm{cheb}(Z_1) = \mathrm{cheb} (\{z,x_2,\ldots,x_N\}$. Hence, we can assume that $x_1 \in S(\mathrm{cheb}(Z_1),R(Z_1))$.
Let us put in correspondence the $N$--net $Z(s) = \{x(s),x_2,\ldots,x_N\}$ and the convex polygon $Q(s)$ defined by its vertices $Z(s)\cap S(\mathrm{cheb}(Z(s)),R(Z(s)))$ to any $s \in [0,x_1y_1]$.
Now for any $s \in [x_1y_1]$ divide the polygon $Q(s)$ into simplices 
$\{\bigtriangleup_1(s),\ldots,\bigtriangleup_{k(s)}(s)\}$
with the common vertex $x(s)$. Now for any 
$s \in [x_1y_1]$ $\mathrm{cheb}(Q(s)) = \mathrm{cheb}(Z(s))$ either belongs to the interior of one simplex or to the common plane of the two of that of the division $\{\bigtriangleup_1(s),\ldots,\bigtriangleup_{k(s)}(s)\}$ (the degenerate cases are also possible).
Consider two cases:

1. Let $\mathrm{cheb}(Z_1)$ belong to the interior of the simplex  
$$
(x_1,a_1,\ldots,a_{N-1}) \in \{\bigtriangleup_1(0),\ldots,\bigtriangleup_{k(0)}(0)\}.
$$
Then---since the Chebyshev center depends on the set is strong continuous way---one can find a minimal
$s_1 \in (0,x_1y_1]$ such that either 
$\mathrm{cheb}(Z_2) = \mathrm{cheb}(\{y_1,a_1,\ldots,a_{N-1}\})$ for $s_1 = x_1y_1$,
or $\mathrm{cheb}(Z(s_1)) = \mathrm{cheb}(\{x(s_1),a_1,\ldots,a_{N-1}\})$ belongs to the common face of one or two simplices of the division 
$\{\bigtriangleup_1(s_1),\ldots,\bigtriangleup_{k(s_1)}(s_1)\}$ if
$s_1 \in (0,x_1y_1)$.

2. Let $\mathrm{cheb}(Z_1)$ belong to the face of the simplex 
$(x_1,a_1,\ldots,a_{N-2},b)$, 
$$
(x_1,a_1,\ldots,a_{N-2},c) \in \{\bigtriangleup_1(0),\ldots,\bigtriangleup_{k(0)}(0)\}.
$$
Then again as in the first case there exists a minimal $s_1 \in (0,x_1y_1]$ such that either
$\mathrm{cheb}(Z_2) = \mathrm{cheb}(\{y_1,a_1,\ldots,a_{N-2},b\})$ or
$\mathrm{cheb}(Z_2) = \mathrm{cheb}(\{y_1,a_1,\ldots,a_{N-2},c\})$ 
if $s_1 = x_1y_1$,
or $\mathrm{cheb}(Z(s_1)) = \mathrm{cheb}(\{x(s_1),a_1,\ldots,a_{N-2},b\})$ or 
$\mathrm{cheb}(Z(s_1)) = \mathrm{cheb}(\{x(s_1),a_1,\ldots,a_{N-2},c\})$ lies on the face of one or two simplices of the division $\{\bigtriangleup_1(s_1),\ldots,\bigtriangleup_{k(s_1)}(s_1)\}$ if 
$s_1 \in (0,x_1y_1)$.
 
In either case lemma 7 and compactness of the interval $[x_1, y_1]$ imply existence of the constant $L_1 >0$ such that
$\mathrm{cheb}(Z_1)\mathrm{cheb}(Z(s_1)) \leq L_1 \alpha (Z_1,Z(s_1))$.
If $s_1 \neq x_1y_1$ then we consider a similar construction for $\mathrm{cheb}(Z(s_1))$ and get the number $s_2$. Continuing the process we get the constants 
$0 < s_1 < \ldots < s_i = x_1y_1$, $L_1,\ldots,L_i > 0$ such that
$$
\mathrm{cheb}(Z_1)\mathrm{cheb}(Z_2) \leq \mathrm{cheb}(Z_1)\mathrm{cheb}(Z(s_1)) + 
\mathrm{cheb}(Z(s_1))\mathrm{cheb}(Z(s_2)) + \ldots +
$$
$$
+ \mathrm{cheb}(Z(s_{i-1}))\mathrm{cheb}(Z_2) \leq 
L_1 \alpha (Z_1,Z(s_1)) + \ldots + L_i \alpha (Z(s_{i-1}),Z_2) 
\leq
$$
$$
\leq (L_1 + \ldots + L_i) \alpha (Z_1,Z_2).
$$
This completes the proof of the statement in this special case.
  
$(ii)$ Consider arbitrary $Z_1 = \{x_1,x_2,\ldots,x_N\}$, 
$Z_2 = \{y_1,y_2,\ldots,y_N\} \in B_{\alpha}(M,\varepsilon)$ and put 
$Z_3 =  \{y_1,x_2,\ldots,x_N\},\ldots,Z_{N+1} = \{y_1,y_2,\ldots,y_{N-1},x_N\}$. 
Again by lemma 6, definition of the Hausdorff metric and triangle inequality together imply existence of the constants 
$L_1 \,\ldots, L_{N-1} >0$ such that
$$
\mathrm{cheb}(Z_1)\mathrm{cheb}(Z_2) \leq \mathrm{cheb}(Z_1)\mathrm{cheb}(Z_3) + \mathrm{cheb}(Z_3)\mathrm{cheb}(Z_4) + \ldots
+ 
$$
$$
+ \mathrm{cheb}(Z_{N+1})\mathrm{cheb}(Z_2) \leq L_1 \alpha (Z_1,Z_3) + L_2 \alpha (Z_3,Z_4) + \ldots
+ 
$$
$$
+ L_{N-1} \alpha (Z_{N+1},Z_2) \leq (L_1 + \ldots + L_{N-1}) \alpha (Z_1,Z_2).
$$
Thus $\mathrm{cheb} : (B_{\alpha}(M,\varepsilon),\alpha) \rightarrow X$ is a Lipschitz map. 
This completes the proof.

{\bf Proof of statement 1.}
Let us denote by $r = R(M)$, $R = R(Z)$, $t = \mathrm{cheb}(M)\mathrm{cheb}(Z) - R - r$ and assume that 
$r  \leq R$. Now the definition of the Hausdorff metric and an inclusion 
$\mathrm{cheb} (Z) \in co(Z)$ imply that
$\sqrt{(R + r + t)^2 + R^2} - r \leq \alpha (M,Z)$ if $N > 3$ and
$\sqrt{(R + t)^2 + R^2} \leq \alpha (M,Z)$ if $N = 3$.
On the other hand
$\mathrm{cheb}(M)\mathrm{cheb}(Z) = R + r + t \leq (1 + \sqrt{5})(\sqrt{(R + r + t)^2 + R^2)} - r)/2$.  
These inequalities put together complete the proof.

  {\bf Proof of statement 2.}
$(i)$ In case both angles $\angle (uwv)$, $\angle (uzv)$ are not acute $\mathrm{cheb}(M)\mathrm{cheb}(Z) = 0$. Assume then that the angle $\angle (uwv)$ is acute one.
Consider a $3$--net $Z_1 = \{u,v,z_1\}$, here the point $z_1$ is constructed rotating the point $z$ over the line $\Pi (u,v)$ by angle equal to one adjacent to that between the halfplanes $\Pi_{+}(u,v,w)$ and $\Pi_{+}(u,v,z)$.
Then the conditions of lemma 9 and geometrical conside\-rations imply that 
$\alpha (M,Z) = \alpha (M,Z_1)$ and $\mathrm{cheb}(M)\mathrm{cheb}(Z) \leq \mathrm{cheb}(M)\mathrm{cheb}(Z_1)$.
Consider two angles $\angle (vuw)$ and $\angle (uvw)$ for the $3$--net $M$. 
Assume now that the angle $\angle (vuw)$ is acute or non-zero and 
$\angle (vuw) \leq \angle (uvw)$.
Then $\mathrm{cheb}(M)\omega (u,v)\leq vw/2 \leq \alpha (M,Z_1)/2$.
Another case (the acute angle $\angle (uzv)$) is considered similarly to the stated one. In the rest of the cases $\omega (u,v)\mathrm{cheb}(Z_1) = 0$.
So using all the inequalities found in this proof we get 
$$
\mathrm{cheb}(M)\mathrm{cheb}(Z) \leq \mathrm{cheb}(M)\mathrm{cheb}(Z_1) \leq \mathrm{cheb}(M)\omega (u,v) +
\omega (u,v)\mathrm{cheb}(Z_1) \leq
$$
$$
\leq \alpha (M,Z).
$$ 
Thus statement $(i)$ of lemma 9 holds true.

$(ii)$ Assume without loss of generality that $[v,q]$ is the interval of the minimal length of the intervals $v,q$, $[v,z]$, $[w,q]$, $[w,z]$ such that its intersection with the interior of the set $co(M)\cup co(Z)$ is empty. Then using the definition of Hausdorff metric and the first part of the statement we find the following inequality:
$$
\mathrm{cheb}(M)\mathrm{cheb}(Z) \leq 
\mathrm{cheb}(M)\mathrm{cheb}(\{u,v,q\}) + \mathrm{cheb}(\{u,v,q\})\mathrm{cheb}(Z) \leq
$$
$$
\leq
\alpha (M,\{u,v,q\}) + \alpha (\{u,v,q\},Z) \leq 2 \alpha (M,Z).
$$
This completes the proof of statement2.

\bigskip
\newpage

{\bf Literature}

\medskip

1. {\bf Garkavi, A.L.} The best possible net and the best possible cross-section of a set in a normed space // Am. Math. Soc., Transl., II. Ser. 39, 111-132 (1964); translation from Izv. Akad. Nauk SSSR, Ser. Mat. 26, 87-106 (1962).

2. {\bf  Kuratowski, C.} Topologie I. Espaces metrisables, espaces complets. 2. ed. revue et augmentee. Monografie Matematyczne. T. XX. Warszawa: Seminarium Matematyczne Uniwersytetu. XI, (1948).

3. {\bf Burago, D.;  Burago, Yu.;  Ivanov, S.} A course in metric geometry.
- Graduate Studies in Mathematics. 33. Providence, RI: American Mathematical Society (AMS). xiv, (2001).  

4. {\bf Belobrov, P.K.} On the Chebyshev point of a system of sets //
Izv. Vyssh. Uchebn. Zaved., Mat. 1966, No.6(55), 18-24 (1966).

5. {\bf Sosov, E.N.} The best net, the best section, and the Chebyshev center of bounded set in infinite-dimensional Lobachevskij space //
Russ. Math. 43, No.9, 39-43 (1999); translation from Izv. Vyssh. Uchebn. Zaved., Mat. 1999, No.9, 42-47 (1999)


6. {\bf Sosov E.~N.} On Hausdorff intrinsic metric // Lobachevskii J.
          of Math. - 2001. - V. 8. - P. 185--189.

7. {\bf Lang U., Pavlovic B., Schroeder V.} Extensions of Lipschitz maps into
Hadamard spaces // Geom. Funct. Anal. - 10, - No. 6. - P. 1527--1553 (2000).

8. {\bf Sosov, E.N.} On Metric space of $2$-nets in the nonpositively curved space // Izv. Vyssh. Uchebn. Zaved., Mat., No.10, 57-60 (1999)

9.  {\bf Shirokov, P.A.} A sketch of the fundamentals of Lobachevskian geometry. Prepared for publication by I.N.Bronshtejn. Translated from the first Russian edition by L.F.Boron. With the assistance of W.D.Bouwsma.
Groningen-The Netherlands: P. Noordhoff Ltd. (1964).

\end{document}